\documentclass[a4paper,12pt,reqno]{amsart}
\newtheorem*{prop}{Proposition}
\newtheorem*{Thm}{Theorem}
\begin{document}

\title{Yangian Construction of the Virasoro Algebra}

\author{S.~Z.~Levendorskii}
\address{Rostov Institute for National Economy\\SU-344798
  Rostov-on-Don\\Russia}

\author{A.~Sudbery}
\address{Department of Mathematics\\University of
  York\\Heslington\\York\\England\\YO1 5DD}
\email{as2@york.ac.uk}

\begin{abstract}
We show that a Yangian construction based on the algebra of an
infinite number of harmonic oscillators (i.e. a vibrating string)
terminates after one step, yielding the Virasoro algebra.
\end{abstract}
\maketitle

In this paper we consider a somewhat more general notion of Yangian than
the usual one (as defined by e.g. Chari and Pressley~\cite[Chapter 12]
{Chari-P}). Let $\mathfrak{g}$
be a Lie algebra, with universal enveloping algebra $U(\mathfrak{g})$,
and suppose $\mathfrak{g} \otimes \mathfrak{g}$ contains an element $\Omega$
which commutes with $\Delta_0 (X)=X \otimes 1 + 1 \otimes X$ in
$U(\mathfrak{g})
\otimes U(\mathfrak{g})$, for all $X \in \mathfrak{g}$. We define a
\textit{Yangian}
of $\mathfrak{g}$ to be a bialgebra $Y(\mathfrak{g})$ containing
$U(\mathfrak{g})$
and another
copy of $\mathfrak{g}$, so that $Y(\mathfrak{g})$ is generated by $\mathfrak{g}
\oplus \theta(\mathfrak{g})$ where $\theta$ is a vector space isomorphism,
with coproducts
\begin{align}
\Delta(X)&=\Delta_0 (X) = X \otimes 1 + 1 \otimes X
  \label{eq:cop1}\\
    \Delta (\theta (X))\, &=\,\theta (X)\otimes 1  \, + \, 1 \otimes
    \theta ( X ) \, + \, [ \Omega , \, X\otimes 1]
\label{eq:cop2}
 \end{align}
where the square brackets denote the commutator.

Usually one imposes relations
\begin{equation}
[X, \,\theta (Y)]=\theta ([X,Y])
\label{eq:usualrel}
\end{equation}
as is justified by the following simple calculation:

\begin{prop}
The coproducts $(\ref{eq:cop1})$ and $(\ref{eq:cop2})$ satisfy
$$
[\Delta (X), \, \Delta (\theta (Y) )]\,=\,\Delta \circ \theta ([X,Y])
\quad
\textrm{modulo}\ R \otimes Y(\mathfrak{g} ) + Y(\mathfrak{g})\otimes R
$$
where $R$ is the subspace of $Y(\mathfrak{g})$ spanned by the relations
$(\ref{eq:usualrel})\ ($i.e. by all elements $[X, \theta(Y)] - \theta([X,Y]))$.
\end{prop}
\begin{proof}
\begin{multline*}
[\Delta (X), \,\Delta (\theta (Y))] \,\,= \,\,[X \otimes 1 + 1 \otimes
X, \,\, \theta (Y) \otimes 1 + 1 \otimes \theta (Y)] \\
+\, [X \otimes 1 + 1 \otimes X,\, [\Omega , \,Y \otimes 1]].
\end{multline*}
The first term on the right-hand side is
\begin{multline*}
[X, \theta (Y)] \otimes 1 \, + \, 1 \otimes [X, \theta (Y)]
\,=\, \theta([X,Y]) \otimes 1 \,+ \, 1 \otimes \theta ([X,Y]) \\
\text{modulo}\quad R\otimes Y(\mathfrak{g}) \, + \, Y(\mathfrak{g}) \otimes
R.
\end{multline*}
The second, using the Jacobi identity and the fact that $\Omega$ commutes
with $X \otimes 1 + 1 \otimes X$, is
$$
[\Omega, \, [X \otimes 1 + 1 \otimes X, \, Y \otimes 1]] = [\Omega ,
\, [X,Y] \otimes 1].
$$
Hence their sum is as stated.
\end{proof}

The commutators $[\theta (X), \theta (Y)]$ are constrained but not completely
determined by the requirement of compatibility with the coproducts
(\ref{eq:cop2}),
and in general they involve new linearly independent elements. Further
linearly independent elements are introduced
by higher commutators, so that the size of the Yangian so defined is
considerably
larger than the polynomial algebra on $\mathfrak{g} \oplus \theta
(\mathfrak{g})$.

The purpose of this note is to point out that when
$\mathfrak{g}=\mathfrak{s}$, the
Lie algebra associated with the oscillations of a string, a slight change
in the above construction leads to a Yangian algebra which closes at the
first level $\theta (\mathfrak{g})$, where it yields the Virasoro algebra.

By the string algebra $\mathfrak{s}$ we mean the algebra of an infinite
number of oscillators whose frequencies are integer multiples of a basic
frequency. Thus we take $\mathfrak{s}$ to have generators $a_m \, (m \in
\mathbb{Z} )$ (raising operators if $m>0$, lowering operators if $m<0$,
central if $m=0$) and $H$ (the Hamiltonian), with commutators
\begin{align}
[a_m,\, a_n ] &= ma_0 \delta_{m+n,0} \label{eq:4}\\
[H,\, a_m ] &= ma_m . \label{eq:5}
\end{align}
This Lie algebra has a Casimir tensor
\begin{equation}
 \Omega = \sum_{m \in \mathbb{Z}} a_m \otimes a_{-m}\ + \ a_0 \otimes
H \, + \, H \otimes a_0 .
\label{eq:6}\end{equation}
We introduce the second copy $\theta(\mathfrak{s})$ with generators $b_m
= \theta (a_m), \ K=\theta (H)$ and coproducts
\begin{align}
\Delta (b_m) &= b_m \otimes 1 \, + \, 1 \otimes b_m \, + \varepsilon
[\Omega,\,a_m\otimes
1] \notag\\
&= b_m \otimes 1 \, + \, 1 \otimes b_m \, + \, \varepsilon m(a_m \otimes
a_0 - a_0 \otimes a_m) \label{eq:7}  \\
\Delta(K)&=K\otimes 1+1\otimes K+\delta [\Omega,\,H\otimes 1]\notag\\
&=K\otimes 1 + 1\otimes K - \varepsilon \Omega \notag
\end{align}
where $\varepsilon \in \mathbb{C}$ is a parameter. The Lie brackets of
the
generators $b_m, \, K$ with each other and with the string generators
$a_m,\, H$ are required to be compatible with the coproducts (\ref{eq:7})
and with the usual Lie algebra coproducts $\Delta _0 (a_m)$ and $\Delta_0(H)$
as in (\ref{eq:cop1}). We note that the only elements of $\mathfrak{s}
\oplus \theta (\mathfrak{s})$ with these Lie algebra-like coproducts are
combinations of $a_m, \, H$ and $b_0$: in Hopf algebra teminology, the
primitive Lie subalgebra of $\mathfrak{s} \oplus \theta(\mathfrak{s})$
is $\mathfrak{s} \oplus \langle b_0 \rangle$.

 The requirement of compatibility gives
\begin{align*}
\Delta ([a_m,\,b_n])&=[\Delta (a_m),\,\Delta(b_n)]\\
 &=[a_m,\,b_n]\otimes 1 \,+\, 1\otimes [a_m,\,b_n],
\end{align*}
so that $[a_m, b_n]$ is Lie algebra-like, and
\begin{align*}
\Delta([H,b_n])&=[\Delta(H),\,\Delta(b_n)]\\
&=[H,b_n]\otimes 1 \, + \, 1 \otimes [H,b_n]\, +\, \varepsilon n^2(a_n
\otimes
a_0 - a_0 \otimes a_n)\\
&=[H,b_n]\otimes 1 \, + \, 1 \otimes [H,b_n]\, +\, n\Delta(b_n) - n(b_n\otimes
1 + 1 \otimes b_n)\end{align*}
so that $[H,b_n]-nb_n$ is Lie algebra-like. The simplest choice is to
take it to vanish:
\begin{equation}   [H,b_n]=nb_n.  \end{equation}
Then $[a_m,b_n]$, which belongs to $\mathfrak{s} \oplus \langle b_0 \rangle$,
must be a multiple of $a_{m+n}$; suppose $[a_m, b_n]=\gamma _{mn} a_{m+n}$.
The Jacobi identity between $a_l, \ a_m$ and $b_n$ gives
$$        l \gamma _{mn}=m \gamma_{ln}            $$
so that $\gamma _{mn}=m\gamma _n$, where the $\gamma _n$ are either all
zero or all non-zero. Taking them to be non-zero, we can redefine the
$b_n$ by dividing by $\gamma _n$ to get
$$
[a_m,b_n]=ma_{m+n}
$$
(but note that this requires the coproduct (\ref{eq:7}) to be modified
by replacing
$\varepsilon$ by $\varepsilon_m = \varepsilon /\gamma _m$).

Now the Jacobi identity between $a_l,\ b_m$ and $b_n$ gives
\begin{align*}
[a_l,\, [b_m,b_n]]&=l(m-n)a_{l+m+n}\\
&=[a_l,\,(m-n)b_{m+n}].
\end{align*}
Put
$$ [b_m,b_n]=(m-n)b_{m+n}\, +\, c_{mn};              $$
then $c_{mn}$ commutes with all $a_l$, and $[H,\, c_{mn}]=(m+n)c_{mn}$,
so if we are to avoid introducing new generators $c_{mn}$ must be a function
of $a_0$ and vanish unless $m+n=0$. The Jacobi identity between $b_l,\
b_m$ and $b_n$ with $l+m+n=0$ gives
$$ (m-n)c_{m+n,-m-n}-(m+2n)c_{m,-m}+(2m+n)c_{n,-n}=0          $$
which has the two independent solutions
$$c_{m,-m}=m \quad \text{and} \quad c_{m,-m}=m^3;          $$
hence the commutator between $b_m$ and $b_n$ must be of the form
\begin{equation}
[b_m,b_n]\, = \, (m-n)b_{m+n}\, + \, \delta_{m+n,0}(m^3F+mG)
\label{Virasoro}\end{equation}
where $F$ and $G$ are functions of $a_0$.

Finally, applying the coproduct (\ref{eq:7}) to this commutator gives
\begin{multline*}
m^3\Delta (F) + m\Delta (G) = m^3(F\otimes 1 + 1 \otimes F) + m(G \otimes
1 + 1 \otimes G)\\
 - \varepsilon _m^2 m^3(a_0 \otimes a_0^2 + a_o^2 \otimes a_0).
\end{multline*}
Write
\begin{align*}
 \Delta ^{\prime} (X) &= \Delta (X) - X \otimes 1 - 1 \otimes X,\\
A &= a_0\otimes a_0^2 + a_0^2 \otimes a_0;
\end{align*}
then
$$
m^3 \Delta ^{\prime} (F) + m \Delta ^{\prime} (G) + m^3 \varepsilon_m^2
A = 0.
$$
We can solve these equations with $m=1$ and $2$ to get $\Delta ^{\prime}
(F)$
and $\Delta ^{\prime} (G)$ as multiples of $A$; then, since $A\neq 0$,
the
remaining equations give $\varepsilon_m$ in terms of $\varepsilon_1$ and
$\varepsilon_2$.
The result is
\begin{align*}
\Delta^{\prime}(F) &= \tfrac{1}{3}(\varepsilon_1^2 - 4 \varepsilon_2^2)A,\\
\Delta^{\prime}(G)&=\tfrac{4}{3}(\varepsilon_2^2-\varepsilon_1^2)A,
\end{align*}
$$ \varepsilon_m^2 \,= \,\varepsilon_1^2 \, + \,
\sum_{r=1}^m\frac{2r-1}{3r^2}(\varepsilon_2^2-\varepsilon_1^2).
$$
Expanding $F$ and $G$ in powers of $a_0$ and using
$$
\Delta(a_0^r)=(a_0\otimes 1+1\otimes a_0)^r=\sum_k\binom{r}{k}a_0^k\otimes
a_0 ^{r-k}, $$
we find
\begin{align*}
F&=\alpha a_0 + \tfrac{1}{9}(\varepsilon_1^2 - 4\varepsilon_2^2)a_0^3,\\
G&=\beta a_0 + \tfrac{4}{9}(\varepsilon_2^2 - \varepsilon_1^2)a_0^3
\end{align*}
where $\alpha$ and $\beta$ are arbitrary.

At this stage we see that we have a bialgebra generated by $a_n,\,b_n$
and $H$, so there is no need to consider the generator $K$. In terms of
our definition, it is a subbialgebra of the Yangian which contains the
Virasoro algebra (\ref{Virasoro}).

In order to exhibit a simple example, we take $\alpha = \beta = 0$ and
$\varepsilon_1 = \varepsilon_2$, so that $\varepsilon_m^2$ is independent
of $m$;
say $\varepsilon_m = \varepsilon$. Then we have a bialgebra generated
by $a_m,\,
b_m$ and $H$ with relations
\begin{align*}
[a_m,a_n]&=ma_0   \delta_{m,-m}\\
[a_m,b_n]&=ma_{m+n}\\
[b_m,b_n]&=(m-n)b_{m+n}- \tfrac{1}{3}\varepsilon^2m^3a_0^3\delta_{m+n,0}\\
[H,\,a_m]&=ma_m\\
[H,\,b_m]&=mb_m
\end{align*}
and coproducts
\begin{align*}
\Delta(a_m)&=a_m\otimes 1 + 1\otimes a_m\\
\Delta(b_m)&=b_m\otimes 1 + 1 \otimes b_m + \varepsilon m(a_m\otimes a_0
- a_0\otimes a_m)\\
\Delta(H) &=h\otimes 1 + 1\otimes H.
\end{align*}

To summarise, our general result is
\begin{Thm}
Let $\mathfrak{s}$ be the Lie algebra with generators $a_m,\,(m \in
\mathbb{Z})$
and  $H$ and Lie brackets $(\ref{eq:4})$ and $(\ref{eq:5})$, and let
$U(\mathfrak{s})$
be the universal enveloping algebra of $\mathfrak{s}$, with the usual
bialgebra structure. Then there is a four-parameter family of bialgebras
containing $U(\mathfrak{s})$ and additional generators $b_m\,(m\in
\mathbb{Z})$,
with coproducts of the form
$$
\Delta(b_m)=b_m\otimes a + 1\otimes b_m + m\varepsilon_m(a_m\otimes a_0
-
a_0\otimes a_m)
$$
and relations including $[H,\,b_m]=mb_m$. In any such bialgebra the
$b_m$ satisfy the Virasoro-like relations $(\ref{Virasoro})$.
\end{Thm}

\end{document}